\documentclass[a4paper,11pt]{amsart} 
\usepackage{amsmath,amssymb} 
\usepackage{mathptmx,courier}
\usepackage[scaled=0.95]{helvet}

\title[For which triangles is Pick's formula almost correct?]
{For which triangles is \\ Pick's formula almost correct?}

\author{Michael Eisermann}
\address{Institut Fourier, Universit\'e Grenoble I, France}
\email{Michael.Eisermann@ujf-grenoble.fr}
\urladdr{www-fourier.ujf-grenoble.fr/~eiserm}

\author{Christoph Lamm}
\address{R{\"u}ckertstra{\ss}e 3, 65187 Wiesbaden, Germany}
\email{Christoph.Lamm@web.de}

\date{first version February 2006; this version compiled \today}

\usepackage[a4paper,
            pdfauthor={Michael Eisermann, Christoph Lamm},
            pdftitle={For which triangles is Pick's formula almost correct?},
            colorlinks=true, 
            anchorcolor=black, 
            linkcolor=black, 
            urlcolor=black,
            citecolor=black, 
            pagecolor=black, 
            filecolor=black, 
            menucolor=black, 
            ]{hyperref}

\theoremstyle{plain}
\newtheorem{theorem}{Theorem}[section]

\theoremstyle{definition}
\newtheorem{definition}[theorem]{Definition}
\newtheorem{conjecture}[theorem]{Conjecture}
\newtheorem{question}[theorem]{Question}
\newtheorem{remark}[theorem]{Remark}
\newtheorem{example}[theorem]{Example}

\newcommand{\N}{\mathbb{N}}
\newcommand{\Z}{\mathbb{Z}}
\newcommand{\R}{\mathbb{R}}
\newcommand{\area}{\operatorname{Area}}
\newcommand{\pick}{\operatorname{Pick}}
\newcommand{\conv}[1]{\operatorname{conv}\bigl\{ #1 \bigr\}}
\begin{document}

\begin{abstract}
  We present an intriguing question about lattice points in triangles 
  where Pick's formula is ``almost correct''.  The question has its origin 
  in knot theory, but its statement is purely combinatorial.  
  After more than 30 years the topological question was recently solved, 
  but the lattice point problem is still open.
\end{abstract}

% \headline{15mm}{\eprintinfo}

\vspace*{-3mm}

\maketitle

%%%%%%%%%%%%%%%%%%%%%%%%%%%%%%%%%%%%%%%%%%%%%%%%%%%%%%%%%%%%%%%%%%%%%%%%%%%%%

\section{Almost Pick's formula}

Let $p,q$ be positive integers and consider the triangle
\[
\textstyle
\Delta = \Delta(p,q) := \conv{ \, \bigl(0,0\bigr), \, \bigl(p,0\bigr), \, 
\bigl( p, \frac{q}{p} \bigr) \, } \subset \R^2 .
\]
We count two types of lattice points in $\Z^2$:
\begin{align*}
  \pick(\Delta) := \;
  & \#\left\{\text{interior lattice points, excluding boundaries}\right\} \\
  \textstyle+\frac{1}{2} 
  & \#\left\{\text{boundary lattice points, excluding vertices}\right\}
  + \textstyle\frac{1}{2} 
\end{align*}
If $\frac{q}{p}$ is an integer, then $\Delta$ is 
a \emph{lattice triangle} and Pick's theorem says that
\[
\area(\Delta) = \pick(\Delta).
\]

This equality will no longer hold in general for $\frac{q}{p} \notin \Z$.
Nevertheless, under favourable circumstances, Pick's formula 
can be \emph{almost correct} in the following sense:

\begin{definition}
  Let $p,q$ be positive integers with $q$ even, so that 
  the area of our triangle $\Delta$ is $\frac{1}{2}q \in \Z$.
  We say that Pick's formula is \emph{almost correct} for $\Delta$
  if $\area(\Delta) = \lfloor \pick(\Delta) \rfloor$, where $\lfloor x \rfloor$ 
  designates the integer part of $x \in \R$.
\end{definition}

Notice that our counting formula defines $\pick(\Delta)$
to be an integer or a half-integer. This means that Pick's
formula is almost correct if and only if $\pick(\Delta)$ equals 
either $\area(\Delta)$ or $\area(\Delta) + \frac{1}{2}$.
Here are two typical examples:

\begin{example}
  For $p=5$, $q=18$ we have $\area(\Delta) = 9$ and $\pick(\Delta) = 9$.
\end{example}

\begin{example}
  For $p=5$, $q=4$ we have $\area(\Delta) = 2$ and $\pick(\Delta) = 2 + \frac{1}{2}$.
\end{example}

Starting with $\Delta$ we can consider 
the magnified triangles $r\Delta$ with $r \in \N$:
\[
\textstyle
r\Delta := \conv{ \, \bigl(0,0\bigr), \, \bigl(pr,0\bigr), \, 
\bigl( pr, \frac{qr}{p} \bigr) \, } \subset \R^2 .
\]
Of course, $p\Delta$ is a lattice triangle.
We can now ask for the stronger condition that Pick's 
formula be almost correct for all $r\Delta$ with $r\in\N$.
Notice that Pick's formula is almost correct for all $r\in\N$
if and only if it is almost correct for all $r = 1,2,\dots,p-1$.

\begin{remark}
  It is a frequently studied question to bound the error 
  between the area and a lattice point count, see for instance 
  the chapter ``A lattice-point problem'' in Hardy \cite{Hardy} 
  and the literature cited there.  Our setting can be seen as 
  the inverse problem: we prescribe very strict error bounds 
  and ask which triangles satisfy them.
\end{remark}

%%%%%%%%%%%%%%%%%%%%%%%%%%%%%%%%%%%%%%%%%%%%%%%%%%%%%%%%%%%%%%%%%%%%%%%%%%%%%

\section{The Casson-Gordon families}

For positive integers $p,q' \in \N$ with $p$ odd and $q'$ even, 
one has $q' = 2kp^2 \pm q$ for some $k \in \N$ and $1 < q < p^2$.  Moreover,
\begin{align*}
  \pick\left(\Delta( p, 2kp^2 + q )\right) 
  & = kp^2 + \pick\left(\Delta(p,q)\right) \quad \text{and} \\
  \pick\left(\Delta( p, 2kp^2 - q )\right) 
  & = kp^2 - \pick\left(\Delta(p,q)\right) + \textstyle\frac{1}{2} .
\end{align*}
This shows that Pick's formula is almost correct for $(p,q')$ 
if and only if it is almost correct for $(p,q)$. 
It is thus natural to restrict attention to $q$ with $1 < q < p^2$.

\begin{theorem}[Casson-Gordon \cite{CassonGordon}] \label{prop:KnownSolutions} 
  Let $p,q \in \Z$ be coprime integers with $1 < q < p^2$, $p$ odd, $q$ even.
  Suppose that $p$ and $q$ satisfy one of the following conditions:
  \begin{enumerate}
  \item $q = np \pm 1$ for some $n \in \N$ with $gcd(n,p)=1$, or
  \item $q = n(p \pm 1)$ for some $n \in \N$ with $n \mid 2p \mp 1$, or
  \item $q = n(p \pm 1)$ for some $n \in \N$ with $n \mid p \pm 1$, $n$ odd, or
  \item $q = n(2p \pm 1)$ for some $n \in \N$ with $(p \mp 1)/n$ odd.
  \end{enumerate}
  Then Pick's formula is almost correct for all triangles 
  $r \Delta(p,q)$ with $r \in \N$, in other words,
  $\area(r\Delta) = \lfloor \pick(r\Delta) \rfloor$ 
  for $\Delta = \Delta(p,q)$ and all $r \in \N$.
\end{theorem}

\begin{remark}
  In the presentation given above the four Casson-Gordon
  families may seem rather complicated at first sight.
  They can be reformulated in a more pleasant and symmetric fashion: 
  each $p^2/q$ has a continued fraction representation of one 
  of the following three types:
  $[a_1,a_2,\dots,a_k,\pm 1,-a_k,\dots,-a_2,-a_1]$ with $a_i > 0$,
  or $[2a,2,2b,-2,-2a,2b]$ or $[2a,2,2b,2a,2,2b]$ with $a,b \ne 0$
  (to obtain all examples we also allow negation and reversal of
  these continued fractions).
  See \cite{Kanenobu}, Theorem 6, for a hint on how to prove this 
  for the first family and use direct calculations for the others.
\end{remark}

\section{Knot-theoretic background}

The only known proof of Theorem \ref{prop:KnownSolutions} 
is intricate and highly indirect, but its story is worth telling.
Since the first version of the present note appeared, in February 2006, 
we have been questioned about the knot-theoretic background, 
and so we feel that we should summarize the proof here and 
give a brief account of its long-winding history.
Even though it is not immediately relevant to 
the combinatorial question towards which we are heading, 
we thus take a detour in order to sketch the argument.
We hope that this will serve to better situate 
the result and motivate the question that ensues.  

\begin{proof}[Topological proof of Theorem \ref{prop:KnownSolutions}]
  The proof is a by-product of a profound topological investigation
  by Casson and Gordon in their seminal work \cite{CassonGordon}.
  They apply the Atiyah-Singer $G$-signature theorem in dimension $4$
  in order to establish necessary conditions for a knot $K \subset \R^3$ 
  to bound a ribbon disk $D \subset \R^3$, $\partial D = K$.
  As a corollary (on page 188 in \cite{CassonGordon}) they show that 
  whenever the two-bridge knot represented by the fraction $q/p^2$ 
  is a ribbon knot,  then Pick's formula is almost correct 
  for all triangles $r \Delta(p,q)$ with $r\in\N$.
  This obstruction allows them to exclude many two-bridge 
  knots, by showing that they cannot bound any ribbon disk.

  On the other hand we have the four families displayed above,
  which have already been stated by Casson and Gordon \cite{CassonGordon},
  alas without proof.  Siebenmann \cite{Siebenmann} proved 
  for two of the Casson-Gordon families that the knot $q/p^2$ 
  is a ribbon knot by explicitly constructing a ribbon disk.
  While pursuing a different approach, Lamm \cite{Lamm} reproved 
  and extended Siebenmann's result by giving a unified construction
  showing that all four Casson-Gordon families yield indeed ribbon knots.
  Together with the fundamental result of Casson and Gordon this implies
  that for the above families Pick's formula must be almost correct,
  as stated in the theorem.
\end{proof}

\begin{remark}
  As a historical note, we mention that
  Siebenmann's contribution \cite{Siebenmann} has not been readily 
  available, and thus the details of the constructive part have been 
  completed in published form only recently in \cite{Lamm} and \cite{Lisca}.  
  The fundamental results of Casson and Gordon \cite{CassonGordon}
  have circulated for more than 10 years only in preprint form.
  Fortunately they have been saved from this fate and preserved 
  for posterity in the book by Guillou and Marin \cite{GuillouMarin}.
\end{remark}

\begin{question}
  The proof via knot theory in dimensions $3$ and $4$ 
  may seem far-fetched for a purely combinatorial statement
  that does not even mention knots nor topology in any way.
  Is there a more direct (combinatorial) proof 
  of Theorem \ref{prop:KnownSolutions}?

  Of course, for a fixed pair $(p,q)$ the theorem can 
  easily be verified by a (computer) count of lattice points.
  It is, however, not obvious how to prove the assertion in general.  
  Is there some more satisfactory (number-theoretic) explanation?
\end{question}

%%%%%%%%%%%%%%%%%%%%%%%%%%%%%%%%%%%%%%%%%%%%%%%%%%%%%%%%%%%%%%%%%%%%%%%%%%%%%

\section{Is the list complete?}

Having set the scene, we now come to the main point
of the present note and formulate the delicate inverse question.
Empirical evidence lets us conjecture that the list stated 
in  Theorem \ref{prop:KnownSolutions} is complete.  
More explicitly this means:

\begin{conjecture} \label{conjecture}
  If $p,q \in \Z$ are coprime integers 
  with $1 < q < p^2$, $p$ odd, $q$ even, and Pick's formula is 
  almost correct for all triangles $r\Delta(p,q)$ with $r \in \N$, 
  then the pair $(p,q)$ belongs to one of the four 
  Casson-Gordon families stated above.
\end{conjecture}

This conjecture is already implicit in the 
article of Casson and Gordon \cite{CassonGordon}, 
who verified it for $p \le 105$ on a computer.  
Although the question has been studied by knot theorists ever since 
the preprint of Casson and Gordon appeared in 1974, the above lattice 
point conjecture is still unsolved after more than 30 years.

\begin{remark}
  The topological problem, sketched above, 
  of classifying two-bridge ribbon knots
  has recently been solved by Lisca \cite{Lisca},
  using an independent topological approach 
  avoiding the combinatorial problem.
  Apart from its own geometric appeal, 
  an affirmative answer to Conjecture \ref{conjecture}
  would have an interesting application in knot theory, 
  as indicated in the preceding proof: 
  it would reprove the result of Lisca, 
  by showing that the Casson-Gordon families exhaust all possibilities.
\end{remark}

\begin{remark}
  We have verified the conjecture for $p < 5000$ using the straightforward 
  counting method.  On an Athlon processor running at 2GHz this took less than 
  $2$ days.  Notice, however, that in its na\"ive form an exhaustive search 
  takes time of order $O(n^5)$ and soon becomes too expensive, so certain 
  optimizations are highly recommended.  \footnote{ If you want to check or further 
  optimize our implementation, you can download it at \linebreak
  \url{http://www-fourier.ujf-grenoble.fr/~eiserm/software.html\#pick}.}
\end{remark}

\begin{remark}
  Following Casson and Gordon \cite[p.\ 187]{CassonGordon},
  in a modified formulation taken from Siebenmann \cite{Siebenmann},
  we write $\sigma(p^2,q,r) := 4\left( \area(r\Delta)-\pick(r\Delta) \right) + 1$
  and have 
  \begin{equation} \label{eq:sumformula}
    \sigma(p^2,q,r) = -\frac{2}{p^2}\sum_{s=1}^{p^2-1} 
    \cot\left(\frac{\pi s}{p^2}\right) \cot\left(\frac{\pi q s}{p^2}\right) \sin^2\left(\frac{\pi q r s}{p}\right).
  \end{equation}
  
  The computation can be sped up with the help of the continued fraction 
  for ${p^2}/{q}$ (``Eisenstein method''):  define the numbers 
  $a_i, q_i > 0$ by $q_0 = p^2$, $q_1 = q$ and $q_{i-1}=a_i q_i + q_{i+1}$. 
  For $x \in \R$ define the function $\{x\}$ as $(\text{fractional part of $x$}) - \frac{1}{2}$ 
  (this is not the standard notation) and the function $((x))$ as $\{x\}$, if $x$ is not an integer and $0$ otherwise.
  According to Siebenmann \cite{Siebenmann} we have (modulo a global sign)
  \begin{equation} \label{eq:contfrac}
    \sigma(p^2,q,r) =  \frac{1}{2} 
      \sum_{i=1}^k (-1)^i 
        a_i \left( 1-4\left\{\textstyle\frac{\displaystyle q_i     r}{\displaystyle p}\right\}^2 \right) 
      - \sum_{i=1}^k (-1)^i 
            \left( 1-4\bigl(\bigl(\textstyle\frac{\displaystyle q_i     r}{\displaystyle p}\bigr)\bigr)
                      \bigl(\bigl(\textstyle\frac{\displaystyle q_{i-1} r}{\displaystyle p}\bigr)\bigr) \right)
  \end{equation}
  
  It would be a welcome complement to the existing literature to elucidate 
  and further develop this ansatz. For instance from \eqref{eq:sumformula} 
  we obtain the symmetry $\sigma(p^2,q,r) = \sigma(p^2,q,p-r)$, 
  for $q q' \equiv 1 \pmod{p^2}$ we have $\sigma(p^2,q,r) = \sigma(p^2,q',q r)$,
  and for $q q' \equiv -1 \pmod{p^2}$ similarly $\sigma(p^2,q,r) = -\sigma(p^2,q',q r)$.
  Equation \eqref{eq:contfrac} allows fast computations and is thus well-suited for empirical explorations.  
  Perhaps it can also provide some hints how to attack Conjecture \ref{conjecture}.
\end{remark}

\begin{remark}
  Considering the average of $\sigma(p^2,q,r)$ over $r=1,\ldots,p-1$, Sikora \cite{Sikora} found
  a relationship with the classical Dedekind sum $s(q,p)$. By Theorem \ref{prop:KnownSolutions}
  if $(p,q)$ belongs to one of the Casson-Gordon families then $\sigma(p^2,q,r)=\pm 1$ for $r=1,\ldots,p-1$.
  In particular $\left|\frac{1}{p-1}\sum_{r=1}^{p-1}\sigma(p^2,q,r)\right|\le 1$. He shows that
  \begin{equation*}
  \sum_{r=1}^{p-1}\sigma(p^2,q,r) = 4\cdot s(q,p)-4 p\cdot s(q,p^2)
  \end{equation*}
  in the following way:
  
  Note that in equation \eqref{eq:sumformula} the variable $r$ occurs only in the $\sin^2$-term. Therefore for summing $\sigma(p^2,q,r)$ over
  $r=1,\ldots,p-1$ we need $\sum_{r=1}^{p-1}\sin^2(\frac{\pi q r s}{p})$ for given $p$, $q$, $s$.
  Because $\sum_{r=1}^{p-1}\sin^2(\frac{\pi q r s}{p})=\frac{p}{2}$ if $p \nmid s$ 
  (and the sum vanishes if $p \mid s$) we obtain:
  \begin{equation*}
  \sum_{r=1}^{p-1}\sigma(p^2,q,r) = 
    -\frac{2}{p^2}\sum_{\substack{ 0 < s < p^2 \\ p \, \nmid \, s }}
    \cot\left(\frac{\pi s}{p^2}\right) \cot\left(\frac{\pi q s}{p^2}\right)\frac{p}{2}
    =-4 p\cdot s(q,p^2) + 4\cdot s(q,p).
  \end{equation*}
  Hence, if $(p,q)$ belongs to one of the Casson-Gordon families then 
  \begin{equation} \label{eq:weakercondition}
  \frac{4}{p-1} \left|s(q,p)-p\cdot s(q,p^2)\right| \le 1.
  \end{equation}
  Which $(p,q)$, with $p$ odd and $q$ even,
  satisfy equation \eqref{eq:weakercondition}? 
  For example, for $p=9$ we find the solutions 
  $q=22,56,68,70$ besides the Casson-Gordon families.
\end{remark}

\begin{remark}
  Writing the (mirrored) triangle $\Delta$ in the form 
  $q x + p^2 y \le p q$ with $x,y \ge 0$, we can apply the lattice point
  counting formula of Beck--Robins \cite[Theorem 2.10]{BeckRobins}. 
  Setting $c_{p,q}:=\frac{1}{4}(1+\frac{1}{p^2}+\frac{1}{q})+\frac{1}{12}(\frac{p^2}{q}+\frac{q}{p^2}+\frac{1}{p^2q})$,
  the number of lattice points in the triangle $t \Delta$ is thus:  
  \[
    L(t)= \frac{1}{2}qt^2+\frac{1}{2}t\left(p+\frac{q}{p}+\frac{1}{p}\right)+c_{p,q}+s_{-tpq}(q,1;p^2)+s_{-tpq}(p^2,1;q),
  \]
  where the last two terms denote the Fourier--Dedekind sums defined in \cite{BeckRobins}.
  (In order to stay as close as possible to the notation 
  in the book we denote the magnifying factor by $t$.)
  Because $q$ divides $-tpq$ we have $s_{-tpq}(p^2,1;q)=s_0(p^2,1;q)$ 
  and we can use reciprocity to obtain
  \[
   s_0(p^2,1;q)=-s_0(q,1;p^2)-c_{p,q}+1.
  \]

  Therefore 
  \[
    L(t)= \frac{1}{2}qt^2+\frac{1}{2}t\left(p+\frac{q}{p}+\frac{1}{p}\right)+s_{-tpq}(q,1;p^2)-s_0(q,1;p^2)+1.
  \]

  In order to compute $\pick(t\Delta)$ we count the lattice points 
  on the catheti ($tp$ and $\lfloor tq/p \rfloor$) and on the hypotenuse
  ($\lfloor t/p \rfloor$ because $q$ and $p^2$ are coprime). 
  Using the notation $\{x\} = x - \lfloor x \rfloor$ the result 
  of subtracting half of the lattice points on the boundary from $L(t)$ 
  (and taking care of the vertices in the way we specified) is:
  \begin{equation} \label{eq:difference}
    \pick(t\Delta)=\frac{1}{2}qt^2+\frac{1}{2}\left\{\frac{t}{p}\right\}
    + \frac{1}{2}\left\{\frac{tq}{p}\right\}+s_{-tpq}(q,1;p^2)-s_0(q,1;p^2). % +\frac{1}{2}.
  \end{equation}
  
  This shows that the $t$-variation of $\area(t\Delta)-\pick(t\Delta)$ depends mostly on 
  the term $s_{-tpq}(q,1;p^2)$, the other terms do not contain $t$ or are small. 
  
  For example, we have for $p=11, q=46$ and $t=1,2$ (illustrated in \cite{Lamm}, page 8):
  $\pick(\Delta) = 23 + \frac{1}{22} + \frac{1}{11} + \frac{6}{11} - \frac{2}{11} = 23.5$ and
  $\pick(2\Delta) = 92 + \frac{1}{11} + \frac{2}{11} - \frac{1}{11} - \frac{2}{11} = 92$.
  
  Formula \eqref{eq:difference} can also be expressed in form 
  of Dedekind--Rademacher sums $r_n(q,p^2)$, 
  see exercise 8.10 in \cite{BeckRobins}.
  Analysing Formula \eqref{eq:difference},
  Beck and Pfeifle \cite{BeckPfeifle} obtained 
  partial results concerning Conjecture \ref{conjecture}.
\end{remark}

\begin{remark}
  Extensions of Conjecture \ref{conjecture} are possible: with the definition 
  \[ I(p,q) := \bigl\{ \sigma(p^2,q,r) \mid r=1,\ldots,p-1 \bigr\} \]
  Conjecture \ref{conjecture} now reads:  $(p,q)$ belongs to one of the Casson-Gordon families
  if and only if $I(p,q) = \{1\}$, $\{-1\}$ or $\{-1,1\}$.
  
  For $I(p,q) = \{-3,-1\}$ we find the following families (with parameter $a>0$):
  \begin{itemize}
  \item $C_1(a)=[2a,-8,-2a,2]$,  
  \item $C_2(a)=[2,2a,-2,2,-2a,-6]$,  
  \item $C_3(a)=[6,2a,-2,2,2a,-2]$ (negative reversed fraction of $C_2$), 
  \item $C_4(a)=[2a,2,-2,2,-2,2,-2,2,-2a-2,2]$, and
  \item $C_5(a)$ $=[2a,2,-2,2,-2,2a+2,-2,2,-2,2]$             
  \end{itemize}
  as well as the sporadic case $[6,-4,-2,2]$. 
  
  The case $I(p,q) = \{-3,-1\}$ and a similar set of families for $I(p,q) = \{1,3\}$ 
  seem to be the only ones besides $\{-1,1\}$ for which exactly two values are attained, 
  meaning that we do not find such families for $I(p,q) = \{3,5\}$, for example.
\end{remark}

\end{document}